\documentclass[a4paper,10pt,twoside]{article}
\pdfoutput=1

\usepackage[english]{babel}
\usepackage[utf8]{inputenc}
\usepackage[T1]{fontenc}
\usepackage{csquotes}

\usepackage{mathtools}
\mathtoolsset{mathic} 
\usepackage{cfr-lm}
\usepackage{dsfont}
\usepackage{microtype}
\usepackage[a4paper]{geometry}

\usepackage{booktabs}
\usepackage{tabu}

\usepackage{lastpage}
\usepackage{fancyhdr}
\usepackage{amsmath,amssymb}
\usepackage{graphicx}
\usepackage{tikz}
\usetikzlibrary{calc}
\usetikzlibrary{cd,positioning,shapes}
\usetikzlibrary{decorations.pathmorphing}
\usetikzlibrary{decorations.markings}
\tikzset{
 closed/.style = {
  decoration = {
   markings,
   mark = at position 0.5 with {
    \node[transform shape, xscale = .8, yscale=.4] {/}; }
  },
  postaction = {decorate}
 },
 open/.style = {
  decoration = {
   markings,
   mark = at position 0.5 with {
    \node[transform shape, scale = .7, yshift=-0.25] {$\circ$}; } },
  postaction = {decorate}
 }
}

\usepackage[inline]{enumitem}
\setlist{noitemsep,nosep,listparindent=\parindent}
\setlist[itemize]{label=\guillemotright}
\setlist[enumerate,1]{ref=\thesubsection.\arabic*}
\setlist[enumerate,2]{label=\alph*.,ref=\theenumi.\alph*}

\usepackage[colorinlistoftodos]{todonotes}
\usepackage[backend=biber,doi=false,url=false,isbn=false,%
safeinputenc,style=alphabetic,citestyle=alphabetic,sorting=nyt]{biblatex}
\bibliography{\jobname.bib}

\usepackage[thmmarks,amsmath]{ntheorem}
\usepackage[colorlinks=true, allcolors=blue]{hyperref}
\usepackage{thmtools}

\numberwithin{equation}{subsection}

\declaretheoremstyle[headformat=swapnumber,headpunct={.\ ---},%
headfont=\normalfont\scshape\lsstyle,bodyfont=\itshape,%
spaceabove=0pt,spacebelow=0pt,%
preheadhook={\bigskip}]{theorem}
\declaretheorem[style=theorem,sibling=subsection]{theorem}
\declaretheorem[style=theorem,sibling=subsection]{proposition}
\declaretheorem[style=theorem,sibling=subsection]{lemma}
\declaretheorem[style=theorem,sibling=subsection]{corollary}
\declaretheorem[style=theorem,sibling=subsection]{conjecture}

\declaretheoremstyle[headformat=swapnumber,headpunct={.\ ---},%
headfont=\normalfont\scshape\lsstyle,bodyfont=\normalfont,%
spaceabove=0pt,spacebelow=0pt,%
preheadhook={\bigskip}]{definition}
\declaretheorem[style=definition,sibling=subsection]{definition}

\declaretheorem[style=definition,sibling=subsection]{example}
\declaretheorem[style=definition,sibling=subsection]{remark}
\declaretheorem[style=definition,sibling=subsection]{notation}

\declaretheoremstyle[headpunct={\!.},headfont=\itshape,bodyfont=\normalfont,%
qed=\ensuremath{\square},spaceabove=0pt,spacebelow=0pt]{proof}
\declaretheoremstyle[headpunct={\!.},headfont=\itshape,bodyfont=\normalfont,%
qed=\ensuremath{\square},spaceabove=0pt,spacebelow=0pt]{nonumberproof}
\declaretheorem[style=proof,numbered=no]{proof}

\declaretheoremstyle[headformat=swapnumber,headpunct={.\ ---},%
headfont=\itshape,bodyfont=\normalfont,qed=\ensuremath{\square},%
spaceabove=0pt,spacebelow=0pt,%
preheadhook={\bigskip}]{nproof}

\renewcommand{\paragraph}[1]{\par\bigskip\refstepcounter{subsection}%
 {\normalfont\normalsize\scshape\noindent\thesubsection%
 \ifthenelse{\equal{#1}{}}%
 {}%
 {\ \textls{#1.}}%
 \ ---}%
}

\usepackage[compress]{cleveref}
\crefformat{subsection}{\S#2#1#3}

\def\ZZ{\mathbb{Z}}
\def\QQ{\mathbb{Q}}
\def\RR{\mathbb{R}}
\def\CC{\mathbb{C}}

\def\QQl{\QQ_\ell}

\def\Gal{\textnormal{Gal}}

\def\old{\textnormal{old}}
\def\new{\textnormal{new}}
\def\tra{\textnormal{tra}}
\def\alg{\textnormal{alg}}

\def\NS{\textnormal{NS}}
\def\AJ{\textnormal{AJ}}

\def\rk{\textnormal{rk}}

\def\Res{\textnormal{Res}}

\def\Hom{\textnormal{Hom}}
\def\End{\textnormal{End}}
\def\Aut{\textnormal{Aut}}

\def\Mat{\textnormal{Mat}}
\def\SO{\textnormal{SO}}
\def\GL{\textnormal{GL}}
\makeatletter
\def\Gmwith[#1]{\mathbb{G}_{\textnormal{m},#1}}
\def\Gmwithout{\mathbb{G}_{\textnormal{m}}}
\def\Gm{\@ifnextchar[{\Gmwith}{\Gmwithout}}
\makeatother

\def\Alb{\textnormal{Alb}}
\def\Albdim{\textnormal{Albdim}}

\def\into{\hookrightarrow}

\def\longto{\longrightarrow}

\def\Spec{\textnormal{Spec}}
\def\Mot{\textnormal{Mot}}
\def\AbMot{\textnormal{AbMot}}
\def\HH{\mathrm{H}_{\textnormal{mot}}}

\def\ChMot{\mathcal{M}_{\textnormal{rat}}}
\def\SmPr{\textnormal{SmPr}}
\def\opp{\textnormal{op}}

\def\QHS{\QQ\textnormal{-HS}}
\def\HB{H}
\def\Hl{H_\ell}
\def\sing{\textnormal{sing}}
\def\et{\textnormal{\'et}}

\def\Ch{h}

\def\GB{G_{\textnormal{MT}}}
\def\Gl{G_\ell}
\def\Glc{\Gl^\circ}

\def\CH{\textnormal{CH}}
\def\cl{\textnormal{cl}}

\title{On the cohomology of surfaces with $p_g = q = 2$
and maximal Albanese dimension}
\author{Johan Commelin \and Matteo Penegini}

\begin{document}
\maketitle

\begin{abstract}
 In this paper we study the cohomology of smooth projective complex surfaces~$S$
 of general type with invariants $p_g = q = 2$ and surjective Albanese morphism.
 We show that on a Hodge-theoretic level, the cohomology is described
 by the cohomology of the Albanese variety and a K3~surface~$X$ that we call the \emph{K3~partner} of~$S$.
 Furthermore, we show that in suitable cases we can geometrically construct
 the K3~partner~$X$ and an algebraic correspondence in $S \times X$ that relates the cohomology of~$S$ and~$X$.
 Finally, we prove the Tate and Mumford--Tate conjectures for those surfaces~$S$
 that lie in connected components of the Gieseker moduli space that contain a product-quotient surface.
\end{abstract}

\section{Introduction} 

\paragraph{}
Let $S$ be a smooth projective complex surface
with invariants $p_g(S) = q(S) = 2$,
and assume that the Albanese morphism $\alpha \colon S \to A$
is surjective.
The results of this paper are inspired by the following two observations:
\begin{enumerate}
 \item The induced map on cohomology
  $\alpha^* \colon H^*(A,\ZZ) \to H^*(S,\ZZ)$ is injective.
  The orthogonal complement
  $H^2_{\new} = H^*(A,\ZZ)^\perp \subset H^*(S,\ZZ)$
  is a Hodge structure of weight~$2$
  with Hodge numbers $(1,n,1)$,
  where $n = h^{1,1}(S) - 4$.
  Such a Hodge structure is said to be of \emph{K3~type}.
 \item Let $S'$ be a smooth projective complex surface
  with invariant $p_g(S') = 1$.
  Then Morrison~\cite{MorrIsog} showed that there
  exists a K3~surface~$X'$ together with
  an isomorphism $\iota' \colon H^2(S',\QQ)^\tra \to H^2(X',\QQ)^\tra$
  that preserves the Hodge structure, the integral structure,
  and the intersection pairing.
  (Here $(\_)^\tra$ denotes the \emph{transcendental} part of a Hodge structure,
  that is, the orthogonal complement of the Hodge classes.)
\end{enumerate}
These observations lead to the following questions.

\paragraph{Question A} 
Let $S$ be as before.

\smallskip
{\narrower\it\noindent
 Does there exist a K3~surface~$X$
 together with an isomorphism
 $\iota \colon (H^2_{\new})^\tra_\QQ \to H^2(X,\QQ)^\tra$
 that preserves
 \begin{itemize*}
  \item[(H)] the Hodge structure, or
  \item[(Z)] the integral structure, or
  \item[(P)] the intersection pairing?
 \end{itemize*}
 \par}
\smallskip

\noindent
We give an affirmative answer to this question in \cref{hodge-k3-partner},
showing that there exists an~$X$ and~$\iota$ that satisfy (H),~(Z), and~(P).
The strategy of the proof is the same as for Morrison's result mentioned above.

\paragraph{Question B} 
\label{questionB}
The Hodge conjecture predicts that if $\iota$ satisfies~(H),
then it is algebraic.

\smallskip
{\narrower\it\noindent
 Do there exist $X$ and~$\iota$ satisfying~(H) as above,
 such that $\iota$ is algebraic?
 \par}
\smallskip

\noindent
In general we are not able to answer this question.
However, an intesting class of examples of the surfaces that we consider
is formed by so-called \emph{product-quotients}:
these are surfaces birational to a surface $(C \times D)/G$,
where $C$ and~$D$ are curves equipped with an action by a finite group~$G$.
In these cases we give an affirmative answer to question~B in \cref{alg-k3-partner}.
The strategy boils down to finding appropriate Prym varieties
in the Jacobians of~$C$ and~$D$,
and taking for~$X$ the associated Kummer variety.

(Note: surfaces~$S$ for which there is a positive answer to question~B
are very much related to the notion of \emph{K3~burgers},
as introduced by Laterveer~\cite{L18}.)

\paragraph{Question C} 
\label{questionC}
Since we are not able to settle question~B in general,
we may aim for something weaker,
sitting in between question~A.(H) and question~B.
We use the notion of \emph{motivated cycles}
introduced by Andr\'e~\cite{AndrMoti}
(see \cref{def-motives} for details).

\smallskip
{\narrower\it\noindent
 Do there exist $X$ and~$\iota$ satisfying~(H) as above,
 such that $\iota$ is motivated?
 \par}
\smallskip

\noindent
Once again, we are not able to give an answer to this question in general.
However, we give sufficient conditions for a positive answer to this question.
For example, we show that to decide question~C one may replace~$S$
with any other surface in the same connected component
of the moduli space of surfaces of general type (\cref{motivated-k3-partner}).
In particular, question~C has a positive answer for every surface~$S$
that lies in the same connected component as a product-quotient surface.

If question~C has a positive answer for the surface~$S$,
then we also deduce that the Tate and Mumford--Tate conjectures hold
for models of~$S$ over finitely generated subfields of~$\CC$ (\cref{T-MTC}).

We summarise these results in the following theorem
(the conjunction of \cref{hodge-k3-partner,alg-k3-partner,motivated-k3-partner,T-MTC}).

\begin{theorem} \label{main-theorem} 
 Let $S$ be a smooth projective complex surface
 with invariants $p_g(S) = q(S) = 2$,
 and assume that the Albanese morphism $\alpha \colon S \to A$
 is surjective.
 \begin{enumerate}
  \item Then there exists a K3~surface~$X$
   and an isomorphism of Hodge structures
   $$\iota \colon (H^2_{\new}(S,\QQ))^\tra \to H^2(X,\QQ)^\tra.$$
  \item If $S$ is a product-quotient surface (\emph{cf}.~\cref{PQsurf})
   with group~$G$,
   then there exist $X$ and~$\iota$ as above,
   and an algebraic cycle on $S \times X$ that induces~$\iota$.
  \item If $S$ is in the same connected component of the Gieseker moduli space as a product-quotient surface,
   then $\iota$ is \emph{motivated} (in the sense of Andr\'e)
   and the Tate and Mumford--Tate conjectures hold for~$S$.
 \end{enumerate}
\end{theorem}

\paragraph{Structure of this text} 
In the next section, \emph{``On the classification of surfaces with $p_g=q=2$''},
we very briefly recall what is known for the surfaces under consideration.
It is important to stress that the classification of these surfaces is not yet complete.
Hence we present the state of the art up to now.
We shall pay particular attention to those surfaces which are product-quotients
(here a classification theorem is available) recalling definitions,
important properties and its associated group theoretical data.
Furthermore, we recall what is known about their moduli space.

In \cref{Hodge-K3} we discuss the existence of Hodge-theoretical K3~partners
for all the surfaces with $p_g=q=2$, following Morrison's theory.
In \cref{hodge-k3-partner} we prove point~1 of \cref{main-theorem}.

In \cref{Alg-K3} we discuss the problem to find a geometric description
where possible of the Hodge-theoretical K3~partners, proving point~2 of \cref{main-theorem}.
Indeed, for those surface which are product-quotients
we are able to find an algebraic K3~partner.

Finally in the last section we see how the results obtained can be used
to prove that the Tate and Mumford--Tate conjectures hold for these surfaces.
As already mentioned here we use the notion of \emph{motivated cycles}
introduced by Andr\'e.
\Cref{T-MTC} proves the last point of \cref{main-theorem}.

\paragraph{Acknowledgements} 
The authors are indebted with Bert van Geemen and Ben Moonen for sharing with them
some of their ideas on this subject. They are also grateful to Matteo Bonfanti,
Fabrizio Catanese, Paola Frediani, Robert Laterveer,
and Jennifer Paulhus for useful discussions and suggestions.

The first author was supported by the
Netherlands Organisation for Scientific Research~(NWO)
under project no.~613.001.207
\emph{(Arithmetic and motivic aspects of the Kuga--Satake construction)}
and by the
Deutsche Forschungs Gemeinschaft~(DFG)
under Graduiertenkolleg~1821
\emph{(Cohomological Methods in Geometry)}.

The second author was partially supported by
MIUR PRIN 2015 ``Geometry of Algebraic Varieties'' and also by GNSAGA of INdAM.

\section{On the classification of surfaces with \texorpdfstring{$p_g = q = 2$}{pg = q = 2}} 
\label{classification}

\paragraph{}
The classification of smooth projective complex surfaces
with invariants $p_g = q = 2$ is not complete,
although there has been much progress in recent years.
We give an overview of the current state of the art.
Let $S$ be a minimal surface of general type with $p_g=q=2$,
so that $\chi(S)=1-q+p_g=1$.
Recall the following classical general inequalities:
\begin{itemize}
 \item $K^2_S \leq 9 \chi(S)$ (Bogomolov-Miyaoka-Yau).
 \item $K^2_S \geq 2 p_g$, if $q>0$, (Debarre).
\end{itemize}
These yield $4 \leq K^2_S \leq 9$ under our assumptions.
Except for the case $K^2_S=9$,
examples of such surfaces have been constructred
for every value of $K^2_S$ in this range.
Since we are dealing with irregular surfaces,
i.e., with $q > 0$, a useful tool to study them is the Albanese map.
The \emph{Albanese variety} of $S$ is
the $q$-dimensional variety
$\Alb(S):=H^0(S,\Omega^1_S)^\vee/H_1(S,\ZZ)$.
By Hodge theory, $\Alb(S)$ is an abelian variety.
For a fixed base point $x_0 \in S$, we define the \emph{Albanese morphism} as
\[
 \alpha = \alpha_{x_0}\colon S \longrightarrow \Alb(S),
 \qquad x \mapsto \int^x_{x_0}.
\]
The dimension of $\alpha(S)$ ($\leq q(S) = 2$) is called the
\emph{Albanese dimension} of $S$ and it is denoted by $\Albdim(S)$.
If $\Albdim(S)=2$, we say that $S$ has \emph{maximal Albanese dimension}.
For surfaces with $q(S)=2$, we have two possibilities:
\begin{enumerate}
 \item $\Albdim(S)=1$ and $\alpha(S)$ is a smooth curve of genus~$2$; or
 \item $\Albdim(S)=2$ and $\alpha$ is a generically finite cover
  of an abelian surface.
\end{enumerate}
The first case is completely understood:
we have a classification theorem, see~\cite{P11}.

The second case is still open.
By \cite[Section 5]{Cat13}
the degree~$\deg(\alpha)$ of the Albanese map is a topological invariant.
In \cref{state_of_the_art}
we summarize the state of the art of the classification using $K^2_S$ and $\deg(\alpha)$ as main invariants.
\begin{table}[!ht]
 \begin{tabular}{cccccccccc}
  \textnumero &$K^2_S$ & $\Albdim$ & $\deg(\alpha)$ & $\#$ & $\dim$ &
  Name & \textsc{mtc} & \textsc{pq/sipc} & Reference \\
  \midrule
  1 & $8$ & $1$ & $-$ & $24$ & $3^{15},4^6,5^2,6$ &
  Isog. Prod. &
  & Yes & \cite{P11}\\
  2 & $8$ & $2$ & $2$ & $2$ & $0^2$ &
  &
  & No & \cite{PRR} \\
  3 & $8$ & $2$ & $ \leq 6$ & $4$ & $3^3,4$ &
  Isog. Prod. &
  \checkmark & Yes & \cite{P11} \\
  4 & $7$ & $2$ & $3$ & $1$ & $3$ &
  &
  \checkmark & Yes & \cite{PiPo17,CF18} \\
  5 & $7$ & $2$ & $2$ & ? & ? &
  &
  & ? & \cite{R15} \\
  6 & $6$ & $2$ & $4$ & $1$ & $4$ &
  &
  \checkmark & Yes & \cite{PP14} \\
  7 & $6$ & $2$ & $2$ & $3$ & $4^2,3$ &
  &
  & No & \cite{PP13a} \\
  8 & $5$ & $2$ & $3$ & $1$ & $4$ &
  Chen--Hacon &
  \checkmark & Yes & \cite{PP13b} \\
  9 & $4$ & $2$ & $2$ & $1$ & $4$ &
  Catanese &
  \checkmark & Yes & \cite{P11,CML02}
 \end{tabular}
 \caption{State of the art of the classification
  of minimal complex algebraic surfaces with invariants $p_g = q = 2$.
  We have indicated, where possible,
  the number of families~($\#$) and their dimensions~($\dim$).
  Moreover, we point out if some members of the family
  are product-quotient surfaces~(\textsc{pq}) or (semi)-isogenous to a product of curves~(\textsc{sipc}), see \cref{PQsurf,semiisog}.
  In the last column,
  we give references to more detailed descriptions of the class.
  Finally, as a consequence of the present paper,
  we put a checkmark in the column~\textsc{mtc}
  if we prove the Tate and Mumford--Tate conjectures for a class.}
 \label{state_of_the_art}
\end{table}

\paragraph{}
In the rest of this section we will describe the examples
$\text{\textnumero} = 3,4,6,8,9$ in more detail.

The examples $\text{\textnumero} = 1,2,5,7$
will not be treated in this paper for the following reasons.
The surfaces in $\text{\textnumero} = 1$ do not have a surjective Albanese map,
therefore this case falls outside the scope of this paper.
For the surfaces in $\text{\textnumero} = 5$, at the time of writing,
nothing is known about their moduli spaces.
Despite the abundance of information about their moduli,
for the surfaces in $\text{\textnumero} = 2$ and~$7$
none of the methods we develop here seem to work;
principally this is due to the small dimensions of their moduli.

\begin{definition} \label{PQsurf} 
 Let $G$ be a finite group acting on two compact Riemann surfaces $C_1$, $C_2$
 of respective genera $g_1, g_2 \geq 2$.
 Consider the diagonal action of $G$ on $C_1 \times C_2$.
 In this situation we say for short:
 the action of $G$ on $C_1 \times C_2$ is \emph{unmixed}.
 By \cite{Cat00} we may assume w.l.o.g.~that $G$ acts faithfully on both factors.

 The minimal resolution $S$ of the singularities of $T = (C_1 \times C_2)/G$,
 is called a \emph{product-quotient surface}.
 If the action of $G$ on the product $C_1 \times C_2$ is free we will speak of
 \emph{surfaces isogenous to a product of unmixed type}.
 In this case $T$ is already smooth.
\end{definition}

\begin{definition} \label{semiisog} 
 Let $C$ be a smooth projective curve,
 and let $G$ be a finite subgroup of $\Aut(C) \times \ZZ/2$.
 Assume that there are elements in $G$
 exchanging the two factors of $C \times C$.
 In this case we say that the action is \emph{mixed}.

 The minimal resolution of singularities
 $S \longrightarrow T = (C \times C)/G$ is called a \emph{mixed surface}.
 A \emph{surface isogenous to a product of mixed type}
 is a mixed surface where $G$ acts freely on $C \times C$.

 We denote by $G_0 \lhd G$ the index two subgroup,
 i.e., the subgroup of elements that do not exchange the factors.
 In general the singularities of $T$ are rather complicated.
 Assuming that $G_0$ acts freely,
 i.e., $(C \times C)/G_0$ is a surface isogenous to a product,
 then $T$ is smooth and we call it a \emph{semi-isogenous mixed surface}.
\end{definition}

\paragraph{} 
The families of surfaces $\text{\textnumero}=1,3$ in \cref{state_of_the_art}
are all the surfaces isogenous to a product of umixed type (see \cite{P11})
and a brief description of those with maximal Albanese dimension (hence \text{\textnumero}=3)
can be read from \cref{PCU_table} below.
There one can find the genus of~$C_1$ and~$C_2$, the group~$G$,
and the number branch points with multiplicity
of the covering $C_i \longrightarrow C_i/G$ for $i=1,2$.
Notice first that the curves $C_i/G$ are elliptic curves for $i=1,2$;
second that for a complete description
one needs a system of generators for the group $G$.
The families of surfaces $\text{\textnumero}=6,8$ and~$9$
contain subloci of product-quotient surfaces.
To briefly describe these particular members
we add information about the singularities of $(C_1 \times C_2)/G$ to \cref{PCU_table}.
\begin{table}[!ht]
 \begin{center}
  \begin{tabular}{
    >{$}c<{$}
   >{$}c<{$}
   >{$}c<{$}
   >{$}c<{$}
   >{$}c<{$}
   >{$}c<{$}
   >{$}c<{$}
   >{$}c<{$}
   >{$}c<{$}
  }
   K_S^2 & g(C_1) & g(C_2) & G & \text{branch} & \text{sing.} & \eta &\dim & $\text{\textnumero}$ \\
   \midrule
   8 & 3 & 3 & V_4 & (2^2), (2^2) & \text{--} & 0 &4 & 3 \\
   8 & 3 & 4 & S_3 & (3), (2^2) & \text{--} & 0 &3 & 3 \\
   8 & 3 & 5 & D_4 & (2), (2^2) & \text{--} & 0 &3 & 3 \\
   \midrule
   6 & 4 & 4 & A_4 & (2), (2) & 2 \times \frac{1}{2}(1,1) & 2 &2 & 6 \\
   5 & 3 & 3 & S_3 & (3), (3) & \frac{1}{3}(1,1) + \frac{1}{3}(1,2) & 3 &2 & 8 \\
   4 & 3 & 3 & Q_8 & (2), (2) & 4 \times \frac{1}{2}(1,1) & 4 & 2 & 9 \\
   4 & 3 & 3 & D_4 & (2), (2) & 4 \times \frac{1}{2}(1,1) & 4 & 2 & 9 \\
   4 & 2 & 2 & C_2 & (2^2), (2^2) & 4 \times \frac{1}{2}(1,1) & 4 & 4 & 9
  \end{tabular}
 \end{center}
 \caption{Classification of
 product-quotient surfaces of unmixed type with $p_g = q = 2$.
 In the column labeled~$\eta$ we give the number of irreducible components
 of the exceptional divisor of the minimal resolution of $(C_1 \times C_2)/G$.
 The last column ($\text{\textnumero}$) mentions the component of the Gieseker moduli space
 described in \cref{state_of_the_art} that contains the family of product-quotient surfaces
 as a subfamily.}
 \label{PCU_table}
\end{table}

\Cref{PCM_table} describes two families of mixed surfaces that are of interest to us.
The surfaces described by the first row are surfaces isogenous to a product of mixed type.
The semi-isogenous mixed surfaces $T = (C \times C)/G$
with $p_g(T) = q(T) = 2$ and $K_T^2 > 0$, form $9$ families~\cite{CF18}.
We will not be able to say anything about the existence of algebraic K3~partners
for these surfaces (\cref{Alg-K3}), but we will study them in \cref{Moti-K3}.
Since the questions studied in \cref{Moti-K3} will turn out to be invariant under deformation
we only need to study one subfamily per component of the Gieseker moduli space.
Of the 9 families mentioned above, 7 are subfamilies of components of the Gieseker moduli space
that also contain families of product-quotient surfaces that we already described
(see \cite{CF18} and~\cite{Pi17}).
The only families for which this is not the case are listed in~\cref{PCM_table};
they are  families in $\text{\textnumero} =3, 4$ described in \cref{state_of_the_art}.

\begin{table}[!ht]
 \begin{center}
  \begin{tabular}{
    >{$}c<{$}
   >{$}c<{$}
   >{$}c<{$}
   >{$}c<{$}
   >{$}c<{$}
   >{$}c<{$}
   >{$}c<{$}
  }
   K_T^2  & g(C) & G  & G^0  & \text{branch} & \dim & $\text{\textnumero}$ \\
   \midrule
   8& 2& C_4& C_2 &   \textit{--} & 3 & 3 \\
   \midrule
   7& 4& C_6 & C_3 &  (2,-4)& 3 & 4 \\
  \end{tabular}
 \end{center}
 \caption{The two families of
 mixed surfaces with $p_g = q = 2$
 that are not subfamilies of components of the Gieseker moduli space
 that also contain families of product-quotient surfaces that we described before.}
 \label{PCM_table}
\end{table}

\paragraph{}
As already remarked, the surfaces in component $\text{\textnumero} = 3$
are all isogenous to a product.
Hence the weak rigidity theorem of Catanese~\cite{Cat00} tells us that
for each family their moduli space consists of
one connected irreducible component in the subspace $\mathcal{M}_{8,2,2}$~($\mathcal{M}_{K^2,p_g,q}$)
of the Gieseker moduli space of surfaces of general type $\mathcal{M}_{8,1}$~($\mathcal{M}_{K^2,\chi}$).
Moreover each member of the family is isogenous to a product.

The property of preserving an isotrivial fibration is no longer true
for the families $\text{\textnumero} = 4,6,8,9$.
Indeed, their moduli space is bigger in some sense.
To be precise let us first analyse the families of surfaces
in \cref{PCU_table} with $K^2_S <8$.
These families form an irreducible sublocus of $\mathcal{M}_{K^2,2,2}$
but they sit inside a bigger connected component.

The connected component with $K^2_S=4$ was studied by~\cite{CML02}.
We have that all the three families in \cref{PCU_table} with $K^2_S=4$
belong to the same connected component of dimension~$4$.
The surfaces with $K^2_S=5$ form an irreducible component of dimension~$2$
sitting inside the connected component
of Chen--Hacon surfaces described in~\cite{PP13b}, which has dimension~$4$.
Finally the family of surfaces with $K^2_S=6$
is a two dimensional irreducible component
inside an irreducible component of $\mathcal{M}_{6,2,2}$ of dimension~$4$,
this component is studied in~\cite{PP14}.

The first entry of \cref{PCM_table} is again isogenous to a product, of mixed type.
Hence its moduli space is $3$-dimensional
and is described by the weak rigidity theorem of Catanese~\cite{Cat00}.

The moduli space of the surfaces relative to the second entry of \cref{PCM_table}
is described by Pignatelli and Polizzi in~\cite{PiPo17}.
In this case the moduli space is a
generically smooth, irreducible, open and normal subset
of the Gieseker moduli space $\mathcal{M}_{7,2,2}$.
For the general member of the family the Albanese surface is simple,
but some specific surfaces admit an irrational fibration over an elliptic curve.

\section{Hodge-theoretic K3~partners} 
\label{Hodge-K3}

\paragraph{} 
Let $S$ be a smooth projective complex surface
with invariants $p_g(S) = q(S) = 2$,
and assume that the Albanese morphism $\alpha \colon S \to A$
is surjective.
Recall question~A from the introduction:

\smallskip
{\narrower\it\noindent
 Does there exist a K3~surface~$X$
 together with an isomorphism
 $\iota \colon (H^2_{\new})^\tra_\QQ \to H^2(X,\QQ)^\tra$
 that preserves
 \begin{itemize*}
  \item[(H)] the Hodge structure, or
  \item[(Z)] the integral structure, or
  \item[(P)] the intersection pairing?
 \end{itemize*}
 \par}
\smallskip

\noindent
In \cref{hodge-k3-partner} we give an affirmative answer to this question.

\paragraph{} 
A \emph{Hodge lattice}~$V$ is a free $\ZZ$-module of finite rank,
endowed with a polarised Hodge structure
such that the polarisation on~$V$ makes $V$ into a lattice%
---a symmetric bilinear form on a free $\ZZ$-module of finite rank.
In particular, the weight (as Hodge structure) of~$V$ is always even.

\begin{notation} 
 There is some risk of confusing notation: if $V$ is a Hodge lattice,
 then $V(k)$ may denote either the $k$-th Tate twist of the Hodge structure on~$V$
 or it may denote the $k$-th twist of the lattice structure on~$V$.
 In this paper we use the notation $V(k)$ only for Tate twists of the Hodge structure.
\end{notation}

\begin{definition} 
 Let $V$ be a Hodge lattice of K3~type
 (that is, the Hodge structure on~$V$ is of K3~type).
 A \emph{K3~partner} of~$V$ is a complex K3~surface~$X$
 together with an isomorphism of Hodge structures
 $\iota \colon V^\tra_\QQ \to H^2(X,\QQ)^\tra$.
 Following the terminology of Morrison (page~181 of~\cite{MorrIsog})
 we say that a K3~partner $(X,\iota)$ is \emph{strict}
 if $\iota$ maps the intersection form on~$V^\tra_\QQ$
 to the intersection form on $H^2(X,\QQ)^\tra$.
 The K3~partner $(X,\iota)$ is \emph{integral}
 if $\iota$ is compatible with an isomorphism of integral Hodge lattices
 $V^\tra \to H^2(X,\ZZ)^\tra$.
\end{definition}

\paragraph{} 
Let $\Lambda$ denote the even unimodular lattice
$E_8\langle-1\rangle^2 \oplus U^3$.
The lattice $\Lambda$ goes by the name \emph{K3~lattice},
since there is an isometry $H^2(X,\ZZ) \cong \Lambda$
for every complex K3~surface~$X$.
For an integer~$d$, let $\Lambda_d$ denote the lattice
$E_8\langle-1\rangle^2 \oplus U^2 \oplus \langle -2d \rangle$.
Observe that $\Lambda_d \into \Lambda$.
The signature of~$\Lambda$ is $(3,19)$,
whereas the signature of~$\Lambda_d$ is $(2,19)$.

Recall that if $M$ is a lattice,
then the pairing on~$M$ defines a natural map $M \into M^\vee$
and the cokernel $M^\vee/M$ is called the \emph{discriminant group}~$A_M$.
The minimal number of generators of~$A_M$ is denoted with~$\ell(A_M)$.

\begin{theorem} 
 Let $L$ and~$M$ be even lattices with signatures
 $(s_+,s_-)$ and $(t_+,t_-)$ respectively.
 Assume that $L$ is unimodular.
 Then there exists a unique primitive embedding
 $M \into L$, if the following conditions hold:
 \begin{enumerate}
  \item $t_+ < s_+$ and $t_- < s_-$; and
  \item $\rk(L) - \rk(M) \ge \ell(A_M) + 2$.
 \end{enumerate}
 \begin{proof}
  This is a slightly weaker form of theorem~1.14.4 of~\cite{NikuLatt}.
 \end{proof}
\end{theorem}

The following corollary is part of an observation by Morrison,
see corollary~2.10 of~\cite{MorrK3Pic}.

\begin{corollary} 
 \label{embed-into-lambda}
 Let $M$ be an even lattice
 with signature $(2,n)$ for some integer $0 \le n \le 8$.
 Then there exists a unique primitive embedding $M \into \Lambda$
 into the K3~lattice introduced above.
\end{corollary}

\paragraph{} 
Let $L$ be a lattice with signature $(s_+,s_-)$,
and assume that $2 \le s_+ \le 3$ and~$s_- \le 19$.
Define $\Omega(L) =
\{x \in \mathbb{P}(L \otimes \CC) \mid (x,x) = 0, (x,\bar x) > 0\}$.
Note that $\Omega(L)$ is an analytic open subset
of the quadric in~$\mathbb{P}(L \otimes \CC)$
defined by the equation $(x,x) = 0$.
If $M \into L$ is an embedding of two such lattices,
then there is a natural holomorphic map $\Omega(M) \to \Omega(L)$.

Observe that there is a natural bijection
\[
 \left\{
  \parbox{6cm}{\strut
   Hodge structures of K3~type on $L$ such that
   for every nonzero $x \in L^{2,0}$ one has
   $(x,x) = 0$, $(x, \bar{x}) > 0$, and~$x \perp L^{1,1}$
  \strut}
 \right\} \longto \Omega(L),
\]
obtained by mapping a Hodge structure on~$L$
to the point~$L^{2,0}$ in~$\Omega(L)$.

\paragraph{} 
We denote with $\mathcal{F}_{2d,\mathbb{K},\CC}^{\textnormal{full}}$
the moduli space of degree~$2d$ primitively polarised K3~surfaces
with full level $\mathbb{K}$-structure
for an admissible subgroup $\mathbb{K} \subset \SO(\Lambda_d)(\hat{\ZZ})$;
see~\cite{Rizov} for details.

\begin{proposition} 
 \label{K3-partner-families}
 Let $B$ be a connected algebraic variety over~$\CC$,
 and let $\mathcal{V}/B$ be a polarised variation of
 $\ZZ$-Hodge structures of weight~$2$ of K3~type.
 Let $b$ be a point of~$B$ and
 assume that $\mathcal{V}_b$ is a lattice with signature~$(2,n)$
 that admits a primitive embedding $\mathcal{V}_b \into \Lambda$.
 Then there is an \'etale open $B^\circ \to B$
 such that $b$ lies in the image of~$B^\circ$,
 a K3~space $f \colon \mathcal{X} \to B^\circ$,
 and a morphism of variations of Hodge structures
 $\mathcal{V}|_{B^\circ} \to \mathrm{R}^2f_*\ZZ$
 that is fibrewise a primitive embedding
 and a Hodge isometry on the transcendental lattices.
 \begin{proof}
  Let $\tilde B \to B$ be a universal cover of~$B$,
  and let $\tilde b$ be a point of $\tilde B$ lying above~$b$.
  Let $L$ denote the lattice underlying the fibre of~$\mathcal{V}_{\tilde b}$.
  Fix a primitive embedding $L \into \Lambda$.
  Under the induced map $\Omega(L) \to \Omega(\Lambda)$
  the Hodge structure on $\mathcal{V}_{\tilde b}$
  maps to a point $x \in \Omega(\Lambda)$.
  By the surjectivity of the period map for K3~surfaces
  there is exists a complex K3~surface~$X$,
  with $H^2(X,\ZZ) \cong \Lambda$,
  such that the Hodge structure on $H^2(X,\ZZ)$ corresponds to~$x$.

  Let $d \in \ZZ$ and $\mathbb{K} \subset \SO(\Lambda_d)(\hat{\ZZ})$
  be such that $X$ admits a primitive polarisation of degree~$2d$
  with full level $\mathbb{K}$-structure.
  This means that $[X] \in \mathcal{F}_{2d,\mathbb{K},\CC}^{\textnormal{full}}$.
  Since $L$ has signature $(2,n)$ we get a primitive embedding $L \into \Lambda_d$.
  This yields a diagram
  \[
   \begin{tikzcd}
    \tilde{B} \dar \rar & \Omega(L) \rar & \Omega(\Lambda_d) \dar \\
    B && \mathbb{K}\backslash\Omega(\Lambda_d)
   \end{tikzcd}
  \]
  The composite map $\tilde{B} \to \mathbb{K}\backslash\Omega(\Lambda_d)$
  factors via a finite cover $B_{\mathbb{K}}$ of $B$.
  \[
   \begin{tikzcd}
    \tilde{B} \dar \rar \drar & \Omega(L) \rar & \Omega(\Lambda_d) \dar \\
    B & B_{\mathbb{K}} \lar{\text{fin \'{e}t}} \rar
    & \mathbb{K}\backslash\Omega(\Lambda_d)
   \end{tikzcd}
  \]
  By a theorem of Borel (see thm~6.4.1 of~\cite{HuybK3}) the map
  $B_{\mathbb{K}} \to \mathbb{K}\backslash\Omega(\Lambda_d)$ is algebraic.
  By proposition~3.2.11 of~\cite{RizovTh},
  there is an open immersion
  $\mathcal{F}_{2d,\mathbb{K},\CC}^{\textnormal{full}} \to
  \mathrm{Sh}_{\mathbb{K}}(\SO(\Lambda_d), \Omega(\Lambda_d))$,
  where the latter is the Shimura variety parameterising polarised Hodge
  structures on $\Lambda_d \otimes \QQ$ with a level $\mathbb{K}$-structure.

  We now have a diagram
  \[
   \begin{tikzcd}
    \tilde{B} \dar \rar \drar & \Omega(L) \rar & \Omega(\Lambda_d) \dar \\
    B & B_{\mathbb{K}} \lar{\text{fin \'{e}t}} \rar
    & \mathbb{K}\backslash\Omega(\Lambda_d) \rar
    & \mathrm{Sh}_{\mathbb{K}}(\SO(\Lambda_d), \Omega(\Lambda_d)) \\
    & B^{\circ} \uar[hook,open] \ar{rr}
    && \mathcal{F}_{2d,\mathbb{K},\CC}^{\textnormal{full}} \uar[hook,open] \\
   \end{tikzcd}
  \]
  where the bottom right rectangle is cartesian.
  By our choice of $d$ and~$\mathbb{K}$,
  we know that $[X] \in \mathcal{F}_{2d,\mathbb{K},\CC}^{\textnormal{full}}$
  and $b \in B$
  are in the image of~$B^\circ$ under the respective maps.
  In particular $B^\circ$ is non-empty.
  Pulling back the universal family of K3~surfaces from
  $\mathcal{F}_{2d,\mathbb{K},\CC}^{\textnormal{full}}$
  to~$B^\circ$
  we end up with a K3~space $f \colon \mathcal{X} \to B^\circ$
  and a morphism of variations of Hodge structures
  $\mathcal{V}|_{B^\circ} \to \mathrm{R}^2f_*\ZZ$.
  By construction it is fibrewise a primitive embedding
  and a Hodge isometry on the transcendental lattices.
 \end{proof}
\end{proposition}

\begin{lemma} 
 \label{h11-k2}
 Let $S$ be a smooth projective complex surface
 with invariants $p_g = q = 2$.
 Let $\alpha \colon S \to A$ be the Albanese morphism,
 and assume that $\alpha$ is surjective.
 Define $H^2_{\new} = (H^*(A,\ZZ)/\mathrm{tors})^\perp \subset H^*(S,\ZZ)/\mathrm{tors}$.
 Then $H^2_\new$ has rank $14 - K_S^2$.
 In particular $(H^2_\new)^\tra$ has signature $(2,n)$ with $n \le 12 - K_S^2$.
 \begin{proof}
  By our assumptions we have $\chi_S = p_g - q + 1 = 1$.
  Noether's formula gives
  $e + K_S^2 = 12 \cdot \chi_S = 12$.
  Observe that $e = \sum_{i=0}^4 (-1)^i b_i$,
  where $b_0 = b_4 = 1$, and $b_1 = b_3 = 2q = 4$.
  Therefore
  \[
   b_2 = 12 - K_S^2 + 2\cdot4 - 2\cdot1 = 20 - K_S^2.
  \]
  Finally, $b_2$ is the sum of the ranks of
  $H^2_\new$ and $H^2(A,\ZZ)$.
  The latter has rank $6$, and we conclude that $H^2_\new$ has rank $14 - K_S^2$.
  Since $(H^2_\new)^\tra$ is a transcendental Hodge structure of K3~type,
  it must have signature $(2,n)$, with $n \le \rk(H^2_\new) - 2$.
 \end{proof}
\end{lemma}

\begin{theorem} 
 \label{hodge-k3-partner}
 Let $S$ be a smooth projective complex surface
 with invariants $p_g = q = 2$,
 and let $\alpha \colon S \to A$ and $H^2_{\new}$
 be as in the preceding~\cref{h11-k2}.
 Then there exists a complex K3~surface~$X$
 and a morphism
 $\iota \colon (H^2_{\new})^\tra_\QQ \to H^2(X,\QQ)^\tra$
 that preserves
 \begin{itemize*}
  \item[(H)] the Hodge structure,
  \item[(Z)] the integral structure, and
  \item[(P)] the intersection pairing.
 \end{itemize*}
 In other words, we have a positive answer to question~A of the introduction.
 \begin{proof}
  Since $4 \le K_S^2 \le 9$, \cref{h11-k2} shows that the lattice $(H^2_{\new})^\tra$
  has signature $(2,n)$ with $n \le 8$.
  Note that $(H^2_\new)^\tra$ is even by the Wu formula:
  for every $v \in H^2(S,\ZZ)$ we have $(v,v) \equiv (v,c_1(S)) \pmod 2$;
  and we have $(H^2_\new)^\tra \subset c_1(S)^\perp$
  since $(H^2_\new)^\tra$ is by definition perpendicular to all Hodge classes.
  The result follows from
  \cref{embed-into-lambda} and \cref{K3-partner-families}.
 \end{proof}
\end{theorem}

\section{Algebraic K3~partners} 
\label{Alg-K3}

Let $S$ be a smooth projective complex surface
with invariants $p_g(S) = q(S) = 2$,
and assume that the Albanese morphism $\alpha \colon S \to A$
is surjective.
In this section we attempt to answer question~B of the introduction:

\smallskip
{\narrower\it\noindent
 Does there exist a K3~surface~$X$
 together with an isomorphism
 $\iota \colon (H^2_{\new})^\tra_\QQ \to H^2(X,\QQ)^\tra$
 that is algebraic?
 \par}
\smallskip

\noindent
Since \cref{hodge-k3-partner} provides an affirmative answer to question~A,
the Hodge conjecture predicts a positive answer to question~B as well.
In \cref{alg-k3-partner} we show that this is indeed the case
for certain surfaces that are product-quotients (\emph{cf}. \cref{PQsurf}).
If the essence of the proof must be captured in one sentence,
it would be the following:
the algebraic correspondence inducing~$\iota$ is built from
the Kummer K3~surface associated with a suitable
$2$-dimensional isogeny factor of the product of the Jacobians
of the curves that are used in the construction
of the product-quotient surface~$S$.

\paragraph{Outline of this section} 
This section is the most technical part of this paper.
It is organised as follows.
First we recall some facts about Chow motives of surfaces.
\Cref{decompositionsecond} describes a natural decomposition of~$\Ch^2(S,\QQ)$,
for product-quotient surfaces~$S$ of unmixed type.
In \cref{alg-k3-partner} we give the general
proof for the existence of an algebraic correspondence inducing~$\iota$.
This proof relies on a case-by-case computation,
for which we refer to a \texttt{MAGMA}-script
of which we tabulate the output.
The final part of this section illustrates this proof
by discussing one of the cases in detail, as an example.

\paragraph{Chow motives of surfaces} 
\label{chowmotsurf}
For an introduction to the theory of Chow motives
we refer to the excellent paper~\cite{Scholl} of~Scholl.
Let $\ChMot$ denote the category of Chow motives over~$\CC$.
We recall that $\ChMot$ is an additive, $\QQ$-linear, pseudoabelian category
(theorem~1.6 of~\cite{Scholl}).
There exists a functor $h \colon \SmPr_{/\CC}^{\opp} \to \ChMot$
from the category of smooth projective varieties over~$\CC$
to the category of Chow motives.
If $X$ is a smooth projective variety, and $G$ is a finite group that acts on~$X$,
then we define $h(X/G) = h(X)^G$.
This leads to a satisfactory theory of motives of quotient varieties,
as is explained in~\cite{motquovar}.

We denote with $\CH^i(X)$ (resp.~$\CH^i(M)$) the $i$-th Chow group
of a smooth projective variety~$X$ (resp.~a motive $M \in \ChMot$).
In general, it is not known whether
the K\"unneth projectors~$\pi_i$ are algebraic,
so it does not (yet) make sense to speak of the summand $h^i(X) \subset h(X)$
for an arbitrary smooth projective variety $X/\CC$.
However, a so-called Chow--K\"unneth decomposition does exist
for curves~\cite{Manin_motive},
for surfaces~\cite{Murre_motsurf},
and for abelian varieties~\cite{DenMur}.
For algebraic surfaces there is in fact the following theorem,
which strengthens the decomposition of the Chow motive.
Statement and proof are copied from theorem~2.2 of~\cite{L18}.

\begin{theorem} \label{ChSurf} 
 Let $S$ be a smooth projective surface over~$\CC$.
 There exists a self-dual Chow--K\"unneth decomposition $\{\pi_i\}$ of~$S$,
 with the further property that there is a splitting
 \[
  \pi_2 = \pi_2^\alg + \pi_2^\tra \quad \in \CH^2(S \times S)
 \]
 in orthogonal idempotents, defining a splitting
 $h^2(S) = h^2(S)^\alg \oplus h^2(S)^\tra$ with Chow groups
 \[
  \CH^i(h^2(S)^\alg) =
 \begin{cases}
  \NS(S) & \text{if $i = 1$,}\\
  0 & \text{otherwise,}
 \end{cases} \qquad
 \text{and} \quad \CH^i(h^2(S)^\tra) =
 \begin{cases}
  \CH^2_{\AJ}(S) & \text{if $i = 2$,}\\
  0 & \text{otherwise.}
 \end{cases}
\]
 Here $\CH^2_{\AJ}(S)$ denotes the kernel of the Abel--Jacobi map.
 \begin{proof}
  The Chow--K\"unneth decomposition is given in proposition~7.2.1 of~\cite{KMP}.
  The further splitting into an algebraic and transcendental component
  is proposition~7.2.3 of~\cite{KMP}.
 \end{proof}
\end{theorem}

\begin{proposition}\label{decompositionsecond} 
 Let $S$ be a product-quotient surface of unmixed type
 with curves $C_1$ and~$C_2$, and with group~$G$.
 Then the second K\"unneth component of~$S$
 is given by $h^2(S)\cong U\oplus Z \oplus E$, where
 \begin{align*}
  U &:= \left(h^2(C_1) \otimes h^0(C_2)\right) \oplus \left(h^0(C_1) \otimes h^2(C_2)\right), \\
  Z &:= \left(h^1(C_1) \otimes h^1(C_2)\right)^G, \\
  E &:= \QQ(-1)^{\oplus \eta},
 \end{align*}
 and where $\eta$ is the number of exceptional divisors introduced in the minimal desingularization of the quotient surface.
 \begin{proof}
  Let $\eta$ denote the number of exceptional divisors introduced in the minimal desingularization of the quotient surface.
  Recall that for the surfaces that we are interested in we gave the value of~$\eta$ in \cref{PCU_table}.
  Observe that $h^2(S) = h^2((C_1 \times C_2)/G) \oplus E$.
  Notice that if the action of $G$ on $C_1\times C_2$ is free then $E=0$.
  By the Künneth formula we obtain
  $$h^2(C_1 \times C_2) = \left(h^2(C_1) \otimes h^0(C_2)\right) \oplus \left(h^0(C_1) \otimes h^2(C_2)\right)
  \oplus \left(h^1(C_1) \otimes h^1(C_2)\right).$$
  By definition we have $h^2((C_1 \times C_2)/G) = h^2(C_1 \times C_2)^G$;
  and since $G$ acts trivially on $h^0(C_i)$ and $h^2(C_i)$, we get the result.
 \end{proof}
\end{proposition}

\begin{remark} \label{decomposefurther} 
 Note that we can further decompose the Chow motive~$Z$ of \cref{decompositionsecond}
 as $Z \cong Z_1 \oplus Z_2$ where
 \[
  Z_1 :=h^1(C_1)^G \otimes h^1(C_2)^G, \quad
  Z_2:=\bigoplus_{W \in \hat{G} -\chi_1} \left(h^1(C_1)^{(W)} \otimes h^1(C_2)^{(W^\vee)}\right)^G,
 \]
 where $\hat{G}$ is the set of isomorphism classes
 of irreducible representations of~$G$ over~$\QQ$,
 $\chi_1$ is the trivial representation,
 and $(\_)^{(W)}$ denotes the $W$-isotypical component.
\end{remark}

\begin{remark} 
 Observe that if $W$ is a finite-dimensional $\QQ$-vector space,
 then we may view it as a Chow motive, as follows:
 It is determined by the identity $\Hom_{\ChMot}(M,W) = \Hom_{\ChMot}(M,\QQ) \otimes W$,
 and it is non-canonically isomorphic to $\QQ^{\oplus \dim(W)}$.

 If $G$ is a finite group,
 and $W$ is equipped with a finite-dimensional representation of~$G$,
 then we may view the representation $G \to \Aut(W)$
 as an action of $G$ on the Chow motive~$W$.
\end{remark}

\paragraph{} 
We now state the main theorem of this section,
which gives a partial answer to question~B of the introduction.
The rest of this section is dedicated to its proof.
We refer to \cref{ex_explicit} for an explicit computation
that illustrates the general argument of this proof.

\begin{theorem} \label{alg-k3-partner} 
 Let $S$ be a minimal surface of general type with $p_g = q = 2$,
 of maximal Albanese dimension isogenous to a product of unmixed type.
 Then there exists a K3~surface~$X$ and a correspondence in $S \times X$ that induces an isomorphism
 between $\HB^{2}(X,\QQ)^\tra$ and $\HB^2(S,\QQ)_\new^\tra$.
\end{theorem}


\paragraph{} 
The proof is done in several steps, and will be completed in \cref{proofcomplete}.
Following \cref{decompositionsecond} and \cref{decomposefurther} we decompose
$h^2(S) = U \oplus Z_1 \oplus Z_2 \oplus E$.
First of all we see that for $i = 1,2$ the quotients $C_i/G=:E_i$ are elliptic curves.
The product $E_1 \times E_2$ is the Albanese variety of~$S$,
and we have $h^2(E_1 \times E_2) = U \oplus Z_1$.
Secondly, let $\eta$ be the number of irreducible components of the exceptional divisor
of the minimal resolution of the quotient surface
(see \cref{PCU_table,PCM_table}; if there is no branch locus, then $\eta = 0$).
This gives $E = \QQ(-1)^\eta$.
For the purpose of this theorem we are interested in the remaining term~$Z_2$.
To find an algebraic K3~partner~$X$ for~$S$,
we will find an abelian surface~$A$ as isogeny factor of~$J(C_1) \times J(C_2)$,
such that $Z_2 = h^2(A)^\tra \oplus \QQ(-1)^{k}$ for some $k \in \ZZ_{\ge 0}$.
We may then take the minimal resolution of singularities
of the Kummer surface $\text{Km}(A) = A/\langle-1\rangle$
for the K3~surface~$X$.
To find the isogeny factor $A$,
we proceed by decomposing (up to isogeny)
the Jacobians $J(C_i)$ of $C_i$ for $i=1,2$ as products of simple abelian varieties following~\cite{PR17}.

\paragraph{} 
Now let $A$ be an abelian variety of dimension $g$ with a faithful action of a finite group $G$.
There is an induced homomorphism of $\QQ$-algebras
\[
 \rho: \QQ[G] \to \End_{\QQ}(A).
\]
Any element $\alpha \in \QQ [G]$ defines an abelian subvariety
\[
 A^{\alpha} := \textrm{Im} (m\rho(\alpha)) \subset A
\]
where $m$ is some positive integer such that $m\rho(\alpha) \in \End(A)$.
This definition does not depend on the chosen integer $m$.

We will now describe the so-called \emph{isotypical decomposition} of the abelian variety~$A$ with group action by~$G$.
Begin with the decomposition of $\QQ[G]$ as a product of simple $\QQ$-algebras $Q_1 \times \cdots \times Q_r$.
The factors $Q_i$ correspond canonically to the rational irreducible
representations $W_i$ of the group $G$, because each one is generated by a
unit element $e_i \in Q_i$ which may be considered as a central idempotent
of $\QQ[G]$.

The corresponding decomposition of $1 \in \QQ[G]$,
$$
1 = e_1 + \cdots + e_r
$$
induces an isogeny, via $\rho$ above,
\begin{equation} \label{eq2.1}
 A^{e_1} \times \cdots \times A^{e_r} \to A
\end{equation}
which is given by addition. Note that the components $A^{e_i}$ are
$G$-stable complex subtori of $A$ with $\Hom_G(A^{e_i},A^{e_j}) =0$ for $i \neq j$.
The decomposition~\eqref{eq2.1} is the isotypical decomposition mentioned above.

\paragraph{} 
The isotypical components $A^{e_i}$ can be decomposed further,
using the decomposition of $Q_i$ into a product of minimal left ideals.
Fix an $i \in \{1, \dots, r\}$, and
let $W_i$ be the irreducible rational representation of~$G$
corresponding to the idempotent~$e_i$.
We will now recall some facts from representation theory;
see \S12.2 of \cite{Serre_LinRepFinGrp} for details.

Write $D_i$ for the simple algebra $\End_G(W_i)$,
and observe that $Q_i = \Mat_{n_i}(D_i^\circ)$
is a matrix algebra of degree $n_i$ over the opposite algebra,
for some $n_i \in \ZZ_{\ge 1}$.
Recall that the Schur index of~$W_i$ is the degree~$m_i$ of $D_i$ over its centre.
If $\chi_i$ is the character of one of the irreducible summands of $W_i \otimes_{\QQ} \CC$,
then $\deg \chi_i = m_i \cdot n_i$.
There is a set of primitive idempotents $\{\pi_{i,1}, \dots ,\pi_{i,n_i}\}$
in $\Mat_{n_i}(D_i^\circ) = Q_i \subset \QQ[G]$ such that
$$e_i = \pi_{i,1} + \cdots + \pi_{i,n_i}.$$
(We warn the reader that the $\pi_{i,j}$ are not $G$-equivariant,
and hence the abelian subvarieties $A^{\pi_{i,j}}$ are not $G$-stable.)
The abelian subvarieties $A^{\pi_{i,j}}$ are mutually isogenous for $j= 1, \dots, n_i$.
Let $B_i$ be any one of these isogenous factors;
we call it a \emph{reduced factor} of~$A^{e_i}$
and $n_i$ the \emph{multiplicity} of the reduced factor:
$B_i^{n_i} \to A^{e_i}$ is an isogeny.
Replacing the factors in \eqref{eq2.1} for every $i = 1, \dots, r$,
we get an isogeny called the {\it group algebra decomposition} of the $G$-abelian variety $A$
\begin{equation} \label{eq2.4}
 B_{1}^{n_1} \times \cdots \times B_{r}^{n_r} \to A.
\end{equation}

Note that, whereas \eqref{eq2.1} is uniquely determined, \eqref{eq2.4} is not.
It depends on the choice of the $\pi_{ij}$ as well as the choice of the $B_i$.
However, the dimension and the isogeny class of the abelian varieties~$B_i$ is independent of choices.

\begin{remark} \label{trivDi}
 If $D_i = \QQ$, then we get
 a $G$-equivariant isomorphism $h^1(A^{e_i}) \cong h^1(B_i) \otimes W_i$
 of Chow motives, where $G$ acts trivially on~$h^1(B_i)$.
\end{remark}

\paragraph{} \label{R:dims} 
While the factors in \eqref{eq2.4} are not necessarily easy to determine,
we may compute their dimension in the case of a Jacobian variety.
Let $C$ be a compact Riemann surface equipped with an action of a finite group~$G$
and consider the induced action of~$G$ on~$J(C)$.
Define $V$ to be the representation of $G$ on $H_1(X,\ZZ) \otimes_{\ZZ} \QQ$.
We use the same notation as at the beginning of this section,
so the quotient~$C/G$ has genus~$g_0$
and the cover $\pi \colon C \to C/G$ has $r$ branch points $\{q_1, \ldots, q_r\}$
where each $q_i$ has corresponding monodromy~$g_i$.
The tuple $(g_1, \ldots, g_r)$ is called the \emph{generating vector} for the action~\cite{B91}.

We now copy equation~2.14 from~\cite{B91}, and explain the notation afterwards:
the Hurwitz character $\chi_V$ associated to $V$ is
\begin{equation}\label{eq:chi}
 \chi_V=2\chi_{1}+2\left(g_{0}-1\right)\rho_{1^*}+\sum\limits_{i=1}^{r}\left(\rho_{1^*}-\rho_{\langle g_{i}\rangle}\right).
\end{equation}
Here $\chi_{1}$ is the trivial character on $G$,
and $\rho_{1^*}$ is the character of the regular representation.
The character~$\rho_H$ is the induced character on~$G$ of the trivial character of the subgroup~$H$.
(When $H=\langle g_i \rangle$, this subgroup is the stabilizer, or isotropy group,
of a point in the fiber of the branch point~$q_i$.)

With this definition of~$\chi_V$ in place, we have the following equality
\begin{equation}\label{eq:dims}
 \dim B_i=\frac{1}{2}\dim_\QQ \pi_{i,j} V=\frac{1}{2} \langle \psi_{i},\chi_{V} \rangle
\end{equation}
where $\psi_{i}$ denotes the character of the $\QQ$-irreducible
representation of $G$ corresponding to $W_i$.
See \cite{Pau08,LR12} for details.

\paragraph{} 
We will now complete the proof of \cref{alg-k3-partner}.
We may calculate the dimension of the $B_i$'s for each class of surfaces in the statement,
either by hand or using the \verb|MAGMA| script as in \cite{PR17}.
The result of this computation is given in \cref{JAC_table}.
\begin{table}[!ht]
 \centering
 {
 \begin{tabular}{
   >{$}c<{$}
  >{$}c<{$}
  >{$}c<{$}
    >{$}c<{$}
  >{$}c<{$}
 }
  G & \textrm{char}(C_1)-\chi_1 & J(C_1)/E_1 &  \textrm{char}(C_2)-\chi_1 & J(C_2)/E_2 \\
  \midrule
 V_4 & [1, 1, 2],\ [1, 1, 4] & E \times L_1 & [1, 1, 3],\ [1, 1, 4] & E \times L_2 \\
 S_3 & [1, 2, 3] & L_1^2 & [1, 1, 2],\ [1, 2, 3] & E \times L_2^2 \\
 D_4 & [1, 2, 5] & L_1^2 & [1, 1, 3],\ [1, 1, 4],\ [1, 2, 5] & E \times E' \times L_2^2 \\
 \midrule
 A_4 & [1, 3, 4] & L_1^3 & [1, 3, 4] & L_2^3 \\
 \midrule
 S_3 & [1, 2, 3] & L_1^2 & [1, 2, 3] & L_2^2 \\
 \midrule
 Q_8 & [2, 1, 5] & A & [2, 1, 5] & A' \\
 D_4 & [1, 2, 5] & L_1^2 & [1, 2, 5] & L_2^2 \\
 C_2 & [1, 1, 2] & L_1 & [1, 1, 2] & L_2 \\
 \end{tabular}
}
 \caption{The group algebra decomposition of the Jacobian varieties~$J(C_i)$.
 The columns $\text{char}(C_i)-\chi_1$ have to be read in the following way:
 each list $[d,n,k]$ represents a non-trivial isotypical factor of $J(C_i)$
 corresponding to a simple rational representation~$W$ of~$G$.
 Here $d$ is the dimension of a reduced factor,
 $n$ is the multiplicity of the factor and
 $k$ is the identifying number according to the \texttt{MAGMA} character table of the group~$G$
 of an irreducible character appearing in $W \otimes_{\QQ} \CC$.}
 \label{JAC_table}
\end{table}

It is a coincidence that in the table all the characters that appear are actually self-dual
and defined over~$\QQ$.
One can check that for each row in the table, there is only one character~$\chi$
that appears in column $\text{char}(C_1)$ such that
the dual character $\chi^\vee$ appears in column $\text{char}(C_2)$.

Let $W$ be the irreducible rational representation of~$G$ that corresponds to~$\chi$.
We will complete the proof by a case distinction.
First assume that $G \ne Q_8$.
In this case, one may check that the Schur index of~$W$ is~$1$,
and in fact $D = \End_G(W) = \QQ$.
In other words, we are in the situtation of \cref{trivDi}.
Let $L_1$ be a reduced factor of~$J(C_1)$ corresponding with~$\chi$,
and denote with~$L_2$ a reduced factor of~$J(C_2)$ that corresponds to $\chi^\vee$.
Thus we have
\[
 h^1(J(C_1))^{(W)} \cong h^1(L_1) \otimes W, \qquad
 h^1(J(C_2))^{(W^\vee)} \cong h^1(L_2) \otimes W^\vee,
\]
as motives with an action of~$G$.
Consequently, we find
\[
 \left(h^1(J(C_1))^{(W)} \otimes h^1(J(C_2))^{(W^\vee)}\right)^G \cong
 h^1(L_1) \otimes h^1(L_2) \otimes (W \otimes W^\vee)^G \cong
 h^1(L_1) \otimes h^1(L_2).
\]
We conclude that $Z_2 \cong h^1(L_1) \otimes h^1(L_2) \cong h^2(L_1 \times L_2)^\tra \oplus \QQ(-1)^{\oplus k}$ for some $k \in \ZZ_{\ge 0}$,
and thus $L_1 \times L_2$ is the $2$-dimensional isogeny factor~$A$ of $J(C_1) \times J(C_2)$ that we are looking for.

\paragraph{The case $G = Q_8$} 
In the case where $G = Q_8$,
we find that $D = \End_G(W) = \mathbb{H}$ which has Schur index~$2$.
In this case we cannot use the methods employed so far
to prove that the Jacobian is isogenous to a product of elliptic curves as in all the other cases.

Nevertheless, in~\cite{FPP} it is proven that in this case the curves $C_1$ and~$C_2$,
which are of genus $g = 3$, admit a bigger automorphism group.
Indeed, their automorphism group is isomorphic to $(C_4 \times C_2) \rtimes C_2$, which readily contains $Q_8$.
More precisely, in~\cite{FPP} it is shown that the curves $C$ of genus $g=3$ and with automorphism $Q_8$ and $(C_4 \times C_2) \rtimes C_2$
give rise to the same subvariety in the moduli space of curves which is the family~(34) of Table~2 in~\cite{fgp}.
Therefore, we can try to decompose the Jacobian of $C$ using this larger group. 

Performing the calculation relative to this larger group we have (in the notation of table \ref{JAC_table}):
\[ 
(C_4 \times C_2) \rtimes C_2 \qquad [1,2,9] \qquad L_1^2 \qquad [1,2,9] \qquad L_2^2.
\] 
To conclude, notice that the nineth character of  $(C_4 \times C_2) \rtimes C_2$ is not self dual,
but we have to restrict it to $Q_8$.
Recalling that $Q_8 \triangleleft (C_4 \times C_2) \rtimes C_2$,
one sees that all the condition of Problem~5.2 on page~65 of~\cite{FH91} are fulfilled.
Hence the restriction of this character to $Q_8$ is
the only two dimensional irreducible representation, which is self dual.
We remark that the $L_i$ must be CM~elliptic curves,
since $\mathbb{H}$ injects into $M_2(\End(L_i)_{\QQ})$.

Now recall that $\mathbb{H} \otimes \mathbb{H} \cong \End(\chi_1 \oplus \chi_i \oplus \chi_j \oplus \chi_k) \cong M_4(\QQ)$.
We conclude that
$h^1(L_1^2) \otimes h^1(L_2^2) \cong
\left( h^1(L_1) \otimes h^1(L_2) \right) \otimes (\chi_1 \oplus \chi_i \oplus \chi_j \oplus \chi_k)$
as Chow motives with an action of~$G$.
In particular we have $\left( h^1(L_1^2) \otimes h^1(L_2^2) \right)^G \cong h^1(L_1) \otimes h^1(L_2)$.
This concludes the proof of \cref{alg-k3-partner}.
\label{proofcomplete}

\paragraph{The case $K^{2} = 8$ and $G = V_4$}\label{ex_explicit} 
As an illustration of the proof of \cref{alg-k3-partner}
and in particular the computations performed by the \texttt{MAGMA} script,
we will now study one example in detail.
This example will occupy us for the next few pages.

Let $E_{1}$ and $E_{2}$ be elliptic curves, marked with two points $\{p_{i,1}, p_{i,2}\} \subset E_{i}$.
As group we take $G=V_4$.
Let $\{\alpha_i, \beta_i\}$ be generators of $\pi_1(E_i,0)$ and let $\gamma_{i,j}$ be a loop in $\pi_{1}(E_{i},0)$ around $p_{i,j}$.
Then we define two $G$-covers $f_{i} \colon C_{i} \to E_{i}$, using the Riemann existence theorem (see e.g., \cite[Sec. III]{Mi95}), by the following epimorphisms of groups:
\begin{equation}\label{}
 \begin{aligned}
  \pi_{1} (E_{1} - \{p_{1,j}\},0) &\longto G \\
  \alpha_1 & \longmapsto (0,1) =: a_1 \\
  \beta_1 & \longmapsto (0,0) =: b_1 \\
  \gamma_{1,j} &\longmapsto (1,0) =: c_{1,j} \\
 \end{aligned}
 \qquad
 \begin{aligned}
  \pi_{1}(E_{2} - \{p_{2,j}\},0) &\longto G \\
  \alpha_2 & \longmapsto (1,0) =: a_2 \\
  \beta_2 & \longmapsto (0,0) =: b_2 \\
  \gamma_{2,j} &\longmapsto (0,1)=: c_{2,j} \\
 \end{aligned}
\end{equation}
By construction $\{p_{i,1}, p_{i,2}\} \subset E_{i}$ is the branch locus, and above each $p_{i,j}$, there are $2$ ramification points with branching orders $2$.
Therefore, by the Riemann--Hurwitz formula the $C_{i}$ are curves of genus $3$.

For a non trivial element $g \in G$, let $\chi_{g}$ denote the nontrivial character of $G$, that annihilates~$g$,
let $\rho_{\langle g \rangle}$ denote the character of $G$ induced from the trivial character of the subgroup generated by~$g$,
let $\rho_{1^*}$ be the regular character.
We get
\[
 \QQ[G] = \QQ^{\chi_{(0,0)}} \times \QQ^{\chi_{(0,1)}} \times \QQ^{\chi_{(1,0)}} \times \QQ^{\chi_{(1,1)}}.
\]

We proceed by calculating the Hurwitz character~\eqref{eq:chi} relative to the first quotient ($i = 1$).
Starting from the ramification data of the curve $C_1$, we get
\[
 \chi_{1}=(1,1,1,1), \quad \rho_{1^*}=(4,0,0,0) \quad \rho_{\langle 1,0 \rangle}=(2,2,0,0),
\]
where the induced trivial representation is calculated using the formula of exercise 3.19.b in~\cite{FH91}.
Therefore, the Hurwitz character is $\chi_V=(6,-2,2,2)$.
Now we use \cref{eq:dims} to compute that the dimensions of the reduced factors~$B_\chi$ of $J(C_1)$
are respectively $1,1,0,1$ for $\chi = \chi_{(0,0)}, \chi_{(0,1)}, \chi_{(1,0)}, \chi_{(1,1)}$.
Analogously, we compute that the dimension of the reduced factors~$B_\chi$ of~$J(C_2)$ are
respectively $1,0,1,1$ for $\chi = \chi_{(0,0)}, \chi_{(0,1)}, \chi_{(1,0)}, \chi_{(1,1)}$.

These results relate to row~1 of \cref{JAC_table} in the following way:
$[1,1,2]$ means that there is a $1$-dimension reduced factor with multiplicity~$1$
corresponding to the second character in the \texttt{MAGMA} character table of~$V_4$,
this is the character~$\chi_{(0,1)}$.
This is exactly the isogeny factor $B_{\chi_{(0,1)}}$ of~$J(C_1)$.
Similarly $[1,1,3]$ corresponds to the isogeny factor $B_{\chi_{(1,0)}}$ of~$J(C_2)$.
Both $J(C_1)$ and $J(C_2)$ have an isogeny factor~$B_{\chi_{(1,1)}}$,
which corresponds to $[1,1,4]$ in \cref{JAC_table},
and thus to the elliptic curves~$L_1$ and $L_2$.



\paragraph{} 
Now we will explain how to construct geometrically the algebraic K3~partner~$X$
of $S = (C_1 \times C_2)/G$, and a correspondence
that induces the isomorphism $\HB^2_\new(S,\QQ)^\tra \to \HB^2(X,\QQ)^\tra$.
The K3~partner~$X$ will turn out to be
the minimal resolution of Kummer surface associated with $L_1 \times L_2$.

Observe that $f_{i}$ factors as in the diagram:
\begin{equation}
 \begin{tikzcd}
  & C_i \dlar[swap,"\phi_{i,(0,1)}"] \dar["\phi_{i,(1,0)}"] \drar["\phi_{i,(1,1)}"] \\
    C_i/\langle (0,1) \rangle \drar[swap,"\psi_{i,(0,1)}"] & C_i/\langle (1,0) \rangle \dar["\psi_{i,(1,0)}"] & C_i/\langle (1,1) \rangle \dlar["\psi_{i,(1,1)}"] \\
    & E_i
 \end{tikzcd}
\end{equation}
Using the Riemann--Hurwitz formula we compute the following genera for the quotient curves:
\[ g( C_1/\langle (0,1)) = 2, \quad g(C_1/\langle (1,0) \rangle)=1, \quad g(C_1/\langle (1,1) \rangle)=2;
\]
\[ g( C_2/\langle (0,1)) = 1, \quad g(C_2/\langle (1,0) \rangle)=2, \quad g(C_2/\langle (1,1) \rangle)=2.
\]
Pushing the preceding diagram through the Jacobian functor, we obtain the diagram:
\[
 \begin{tikzcd}
  & J(C_i) \dlar \dar \drar \\
  J(C_i/\langle (0,1) \rangle) \drar & J(C_i/\langle (1,0) \rangle) \dar & J(C_i/\langle (1,1) \rangle) \dlar \\
  & E_i
 \end{tikzcd}
\]
This leads to the following isogenies of abelian varieties:
\begin{equation}\label{eq_IntermPryms}
 \begin{aligned}
  J(C_{1}/\langle (0,1)\rangle &\sim E_{1} \times P(\psi_{1,(0,1)}) & J(C_{2}/\langle (0,1)\rangle &\sim E_{2} \\
  J(C_{1}/\langle (1,0)\rangle &\sim E_{1} & J(C_{2}/\langle (1,0)\rangle &\sim E_{2} \times P(\psi_{2,(1,0)}) \\
  J(C_{1}/\langle (1,1)\rangle &\sim E_{1} \times P(\psi_{1,(1,1)}) & J(C_{2}/\langle (1,1)\rangle &\sim E_{2} \times P(\psi_{2,(1,1)}), \\
 \end{aligned}
\end{equation}
where $P(\psi)$ denotes the \emph{Prym--Tyurin variety} associated to the cover $\psi$:
it is the kernel of the induced map $J(\psi)$ between the Jacobians
(see also, \cite{BL} paragraph 12.2).
Observe that $L_i$ is isogenous to $P(\psi_{i,(1,1)})$.
Finally by \eqref{eq_IntermPryms} we have
\begin{equation}
 \begin{aligned}
  J(C_1) & \sim J(C_{1}/\langle (0,1)\rangle \times P(\phi_{1,(0,1)}) \\
  & \sim E_{1} \times P(\psi_{1,(0,1)}) \times P(\phi_{1,(0,1)})
 \end{aligned}
\end{equation}
\begin{equation}
 \begin{aligned}
  J(C_2) & \sim J(C_{2}/\langle (1,0)\rangle \times P(\phi_{2,(1,0)}) \\
  & \sim E_{2} \times P(\psi_{2,(1,0)}) \times P(\phi_{2,(1,0)}).
 \end{aligned}
\end{equation}




We now go back to the surface $S=(C_1 \times C_2)/V_4$.
Since $\eta = 0$, we get
\[
 h^2(S) \cong h^2(E_1 \times E_2) \oplus h^1(P(\psi_{1,(1,1)})) \otimes h^1(P(\psi_{2,(1,1)})).
\]
Let us go further and build an algebraic K3~partner of $S$.
To do that we consider the abelian surface $A = P(\psi_{1,(1,1)}) \times P(\psi_{2,(1,1)})$ and divide modulo the natural involution.
In this way we get a singular Kummer surface.


According to the propostion above and Shioda and Inose~\cite{SI77}
the minimal resolution $X$ of the singularities of $\textrm{Km}(A)$
is a K3~surface whose transcendental part of $H^2(X)$ is isomorphic to the transcendental part of~$H^2(A)$.
Now consider the following diagram:
\[
 \begin{tikzcd}
  C_1 \times C_2 \dar \rar & J(C_1) \times J(C_2) & A \lar \dar & X \dlar \\
  S & & \textrm{Km}(A)
 \end{tikzcd}
\]
All the morphisms in this diagram induce correspondences
and by composing these correspondences we obtain an
isomorphism $h^2(S)^\tra_\new \to h^2(X)^\tra$
that induces an isomorphism of Hodge structures
$\HB^2_\new(S,\QQ)^\tra \to \HB^2(X,\QQ)^\tra$.

\section{Motivated K3~partners} 
\label{Moti-K3}

Let $S$ be a smooth projective complex surface
with invariants $p_g(S) = q(S) = 2$,
and assume that the Albanese morphism $\alpha \colon S \to A$
is surjective.
In this section we attempt to answer question~C of the introduction:

\smallskip
{\narrower\it\noindent
 Does there exist a K3~surface~$X$
 together with an isomorphism
 $\iota \colon (H^2_{\new})^\tra_\QQ \to H^2(X,\QQ)^\tra$
 that is motivated in the sense of Yves Andr\'e?
 \par}
\smallskip

\paragraph{} 
In this section we prove the Tate and Mumford--Tate conjectures
for surfaces that fall into type $\text{\textnumero} = 3, 4, 6, 8$, and~$9$
of \cref{state_of_the_art}.
We will use the language of motives,
and specifically motivated cycles as introduced by Andr\'e~\cite{AndrMoti}.

This section is organised as follows:
First we introduce notation
and recall the definition of the Mumford--Tate group
and the $\ell$-adic monodromy groups.
Then we will recall three conjectures that are connected in the following sense:
if two of the conjectures hold, then so does the third.
These conjectures are
\begin{enumerate*}[label=(\textit{\roman*})]
 \item the Hodge conjecture;
 \item its $\ell$-adic analogue known as the Tate conjecture; and
 \item the Mumford--Tate conjecture.
\end{enumerate*}

Starting from \cref{def-motives}
we recall the definition of motivated cycles in the sense of
Andr\'e~\cite{AndrMoti},
and we quote the main theorems that describe the resulting category of motives.
Once we we have all the machinery in place
we turn our attention to the proof of the Tate and Mumford--Tate conjectures
for the surfaces mentioned above.

\paragraph{Notation} 
Let $K \subset \CC$ be a field,
and let $X$ be a smooth projective variety over~$K$.
We denote with $\HB^i(X)$ the singular cohomology group
$H_{\sing}^i(X(\CC),\QQ)$.
It is naturally endowed with a pure Hodge structure of weight~$i$.
Let $\ell$ be a prime number,
and let $\bar K \subset \CC$ be the algebraic closure of $K \in \CC$.
We denote with $\Hl^i(X)$ the $\ell$-adic \'etale cohomology group
$H_{\et}^i(X_{\bar K},\QQl)$.
It is naturally endowed with a Galois representation
$\Gal(\bar K/K) \to \GL(\Hl^i(X))$.

Artin's comparison theorem between
\'etale cohomology and singular cohomology
gives an isomorphism of vector spaces
\[
 \Hl^i(X) \cong \Hl^i(X_\CC) \stackrel{\sim}{\longrightarrow}
 \HB^i(X) \otimes \QQl
\]
that is functorial in~$X$.

Recall from \cref{chowmotsurf} that
we denote with $\CH^*(X)$ the Chow ring of~$X$ with $\QQ$-coefficients.
Recall the cycle class map
$\cl \colon \CH^i(X_\CC) \to \HB^{2i}(X)(i)$
for singular cohomology.
There is also a cycle class map
$\cl_\ell \colon \CH^i(X) \to \Hl^{2i}(X)(i)$
for \'etale cohomology.
These are compatible with the comparison isomorphism
$\Hl^i(X) \cong \HB^i(X) \otimes \QQl$;
we get the following commutative diagram:
\[
 \begin{tikzcd}
  \CH^i(X) \ar[rr,"\cl_\ell"] \ar[d] && \Hl^{2i}(X)(i) \ar[d,"\sim"] \\
  \CH^i(X_\CC) \ar[r,"\cl"] & \HB^{2i}(X)(i)
  \ar[r,"\_ \otimes 1"] & \HB^{2i}(X)(i) \otimes \QQl
 \end{tikzcd}
\]
\paragraph{Mumford--Tate groups} 
Let $V$ be a $\QQ$-Hodge structure.
There is a representation of $\mathbb{S} = \Res_{\CC/\RR} \Gm$
on $V \otimes \RR$:
on complex points $(z,\bar z)$ acts on $v \in V^{pq}$
by $v \mapsto z^{-p}\bar z^{-q}v$.
(The minus signs are a historical convention.)
Write $h = h_V$ for this representation $\mathbb{S} \to \GL(V)_\RR$.

The \emph{Mumford--Tate group} of~$V$
is the smallest algebraic subgroup $G \subset \GL(V)$ over~$\QQ$
such that $G_\RR$ contains the image of~$h_V$.
We denote the Mumford--Tate group of~$V$ with~$\GB(V)$.
Alternatively, $\GB(V)$ may be defined using the Tannakian formalism.
It is the algebraic group over~$\QQ$
associated with the Tannakian subcategory of $\QHS$ generated by~$V$.
If $V$ is polarisable,
then this subcategory generated by~$V$ is semisimple,
and hence $\GB(V)$ is a reductive algebraic group.

Two more remarks are in place:
First, observe that $\GB(V)$ is a connected algebraic group,
since $\mathbb{S}$ is connected.
Second, note that
the subspace $V \cap V_\CC^{0,0}$ of \emph{Hodge classes} in~$V$
is exactly the space of invariants~$V^{\GB(V)}$.

\paragraph{$\ell$-adic monodromy groups} 
Let $K \subset \CC$ be a field of finite transcendence degree over~$\QQ$.
Let $\ell$ be a prime number,
let $V$ be a finite-dimensional $\QQl$-vector space,
and let $\rho \colon \Gal(\bar K/K) \to \GL(V)$
be a representation that is continuous for the $\ell$-adic topology on~$V$.

The \emph{$\ell$-adic monodromy group} of~$V$
is the smallest algebraic subgroup $G \subset \GL(V)$ over~$\QQl$
such that $G(\QQl)$ contains the image of~$\rho$.
We denote the $\ell$-adic monodromy group of~$V$ with~$\Gl(V)$.
In general, the algebraic group $\Gl(V)$ is not connected;
the identity component is denoted $\Glc(V)$.

An element of $V$ is called a \emph{Tate class}
if it is invariant under an open subgroup of $\Gal(\bar K/K)$.
In particular, the subspace of Tate classes in~$V$
is exactly the space of invariants~$V^{\Glc(V)}$.

In general, the algebraic group $\Glc(V)$ is not reductive.

\begin{conjecture}[Hodge] 
 Let $X$ be a smooth projective variety over~$\CC$.
 Then the image of $\cl \colon \CH^i(X) \to \HB^{2i}(X)(i)$
 is the subspace of Hodge classes
 $\HB^{2i}(X)(i) \cap \HB^{2i}(X)(i)_\CC^{0,0}$.
\end{conjecture}

\begin{conjecture}[Tate] 
 Let $X$ be a smooth projective variety over a number field~$K$.
 Then the image of
 $\cl_\ell \colon \CH^i(X_{\bar K}) \to \Hl^{2i}(X)(i)$
 spans the space of Tate classes $\Hl^{2i}(X)(i)^{\Glc(X)}$.
\end{conjecture}

\begin{conjecture}[Mumford--Tate] 
 Let $X$ be a smooth projective variety over a number field~$K$.
 The comparison isomorphism $\Hl^i(X) \to \HB^i(X) \otimes \QQl$
 induces and isomorphism
 $\Glc(\Hl^i(X)) \to \GB(\HB^i(X)) \otimes \QQl$.
\end{conjecture}

\begin{remark} 
 To illustrate how these conjectures fit together,
 we make the following claims.
 \begin{enumerate}
  \item If the Mumford--Tate conjecture is true for~$X$,
   then the Hodge conjecture for~$X$
   is equivalent to the Tate conjecture for~$X$.
  \item If the Tate conjecture is true for all
   smooth projective varieties~$X$ over~$K$,
   then the $\ell$-adic monodromy groups are reductive.
   This follows from \cite{Moonen_remarkTate}.
  \item If the Hodge and Tate conjectures are true for all~$X$,
   then the Mumford--Tate conjecture is true for all~$X$.
 \end{enumerate}
\end{remark}

\paragraph{Motivated cycles} 
\label{def-motives}
Let $K$ be a subfield of~$\CC$,
and let $X$ be a smooth projective variety over~$K$.
A class $\gamma$ in $H^{2i}(X)$ is called
a \emph{motivated} cycle of degree~$i$
if there exists an auxiliary smooth projective variety~$Y$ over~$K$
such that $\gamma$ is of the form $\pi_*(\alpha \cup \star\beta)$,
where $\pi \colon X \times Y \to X$ is the projection,
$\alpha$ and~$\beta$ are algebraic cycle classes in $H^*(X \times Y)$,
and $\star\beta$ is the image of~$\beta$ under the Hodge star operation.
(Alternatively, one may use the Lefschetz star operation,
see~\S1 of~\cite{AndrMoti}.)

Every algebraic cycle is motivated,
and under the Lefschetz standard conjecture the converse holds as well.
The set of motivated cycles naturally forms a graded $\QQ$-algebra.
The category of motives over~$K$, denoted~$\Mot_K$,
consists of objects $(X,p,m)$,
where $X$ is a smooth projective variety over~$K$,
$p$ is an idempotent motivated cycle on $X \times X$,
and $m$ is an integer.
A morphism $(X,p,m) \to (Y,q,n)$
is a motivated cycle~$\gamma$ of degree $n-m$ on $Y \times X$
such that $q \gamma p = \gamma$.
We denote with $\HH(X)$ the object $(X,\Delta,0)$,
where $\Delta$ is the class of the diagonal in $X \times X$.
The K\"unneth projectors $\pi_i$ are motivated cycles,
and we denote with $\HH^i(X)$ the object $(X,\pi_i,0)$.
Observe that $\HH(X) = \bigoplus_i \HH^i(X)$.
This gives contravariant functors $\HH(\_)$ and $\HH^i(\_)$
from the category of smooth projective varieties over~$K$ to~$\Mot_K$.

\begin{theorem} 
 \label{mot-props}
 The category $\Mot_K$ is Tannakian over~$\QQ$,
 semisimple, graded, and polarised.
 Every classical cohomology theory of smooth projective varieties over~$K$
 factors via~$\Mot_K$.
 \begin{proof}
  See th\'eor\`eme~0.4 of~\cite{AndrMoti}.
 \end{proof}
\end{theorem}

\begin{definition} 
 Let $K$ be a subfield of~$\CC$.
 An \emph{abelian motive} over~$K$ is an object
 of the Tannakian subcategory of~$\Mot_K$
 generated by objects of the form $\HH(X)$
 where $X$ is an abelian variety,
 or $X = \Spec(L)$ for some finite extension
 $L/K$, with $L \subset \CC$.

 We denote the category of abelian motives over~$K$ with~$\AbMot_K$.
\end{definition}

\begin{example} 
 \label{K3-abelian-motive}
 If $X/K$ is a K3~surface,
 then $\HH^2(X)$ is an abelian motive,
 by th\'eor\`eme~7.1 of~\cite{AndrMoti}.

 The Lefschetz motive $\QQ(-1)$ is abelian,
 because any class of a hyperplane section in an abelian variety~$A$
 will give a splitting $\HH^2(A) \cong M \oplus \QQ(-1)$.
\end{example}

\begin{theorem} 
 \label{hodge-is-motivated}
 The Hodge realisation functor $\HB \colon \AbMot_\CC \to \QHS$
 is a full functor.
 \begin{proof}
  See th\'eor\`eme~0.6.2 of~\cite{AndrMoti}.
 \end{proof}
\end{theorem}

\begin{theorem} 
 \label{variational-hodge-is-motivated}
 Let $B$ be a reduced connected scheme of finite type over~$\CC$.
 Let $f \colon X \to B$ be a smooth projective morphism,
 and $\xi$ a global section of the sheaf $\mathrm{R}^{2i}f_*\QQ(i)$.
 If there is a point $0 \in B(\CC)$ such that $\xi_0$ is motivated,
 then $\xi_b$ is motivated for all $b \in B(\CC)$.
 \begin{proof}
  See th\'eor\`eme~0.5 of~\cite{AndrMoti}.
 \end{proof}
\end{theorem}

\paragraph{} 
By \cref{mot-props}, the singular cohomology and $\ell$-adic cohomology functors
factor via~$\Mot_K$.
This means that if $M$ is a motive,
then we can attach to it a Hodge structure~$\HB(M)$
and an $\ell$-adic Galois representation $\Hl(M)$.
The comparison isomorphism between singular cohomology
and $\ell$-adic cohomology extends to an isomorphism of vector spaces
$\Hl(M) \cong \HB(M) \otimes \QQl$ that is natural in the motive~$M$.

We shall write $\GB(M)$ for $\GB(\HB(M))$.
Similarly, we write $\Gl(M)$ (resp.~$\Glc(M)$)
for $\Gl(\Hl(M))$ (resp.~$(\Glc(\Hl(M)$).
The Mumford--Tate conjecture extends to motives:
for the motive~$M$ it asserts that the comparison isomorphism
$\Hl(M) \cong \HB(M) \otimes \QQl$
induces an isomorphism
$\Glc(M) \cong \GB(M) \otimes \QQl$.

\paragraph{} 
We now have the notation and theory in place to answer question~C
(see~\cref{questionC})
about surfaces of general type with $p_g = q = 2$.
We give a partial answer to this question in the following results.

\begin{theorem} 
 \label{motivated-k3-partner}
 Let $\alpha \colon S \to B$ be a smooth projective family
 of surfaces of general type
 with invariants $p_g = q = 2$ and dominant Albanese morphism.
 Assume that $B$ is connected,
 and assume that there is one point $0 \in B$
 such that the motive $\HH^2(S_0)$ of the fibre~$S_0$ is an abelian motive.
 Then for every point $b \in B$,
 there exists a K3~surface $X_b$,
 and an isomorphism of motives $\HH^2(S_b)_\new^\tra \cong \HH^2(X_b)^\tra$.
 In particular, the motive $\HH^2(S_b)$ is abelian.
 \begin{proof}
  The main idea of the proof is as follows:
  Using \cref{K3-partner-families}
  we construct a family $X \to B$ of Hodge-theoretic K3~partners.
  We then use \cref{hodge-is-motivated}
  to prove that $\HH^2(S_0)_\new^\tra$ is isomorphic to $\HH^2(X_0)_\new^\tra$.
  Finally, this isomorphism spreads out to the other fibres
  via \cref{variational-hodge-is-motivated}.
  We now make this sketch precise.

  By replacing $\alpha \colon S \to B$
  with the pullback along $\tilde B \to B$,
  we may and do assume that $B$ is simply connected.
  Let $\mathcal{V}$ denote the
  subvariation of Hodge structures of~$\mathrm{R}^2\alpha_*\QQ$
  whose fibre $\mathcal{V}_b$ at a point $b \in B$ is $\HB^2(S_b,\QQ)_\new$.
  Fix a point $b \in B$.
  By \cref{K3-partner-families}
  we find that there is an open $B_b^\circ \subset B$
  such that $b \in B_b^\circ$,
  a K3~space $f \colon \mathcal{X} \to B_b^\circ$,
  and a morphism of variations of Hodge structures
  $\iota \colon \mathcal{V}|_{B_b^\circ} \to \mathrm{R}^2f_*\QQ$
  that is fibrewise a primitive embedding
  and a Hodge isometry on the transcendental lattices.
  We may view $\iota$ as a global section
  of the sheaf $\mathcal{V}^\vee|_{B_b^\circ} \otimes \mathrm{R}^2f_*\QQ$
  which is a subsheaf of $\mathrm{R}^4(\alpha,f)_*\QQ(2)$.
  Note that we may and do assume that $b \in B_0^\circ$;
  indeed, if $b \notin B_0^\circ$, then we first prove the statement
  for all points in $B_0^\circ$,
  and then rerun the proof with a point $0' \in B_b^\circ \cap B_0^\circ$.

  Recall from \cref{K3-abelian-motive} that $\HH^2(X_0)$ is an abelian motive.
  Also note that $\HH^2(S_0)$ is abelian by assumption.
  Hence $\iota_0$ is motivated, by \cref{hodge-is-motivated}.
  By \cref{variational-hodge-is-motivated},
  we see that $\iota_b$ is motivated as well.
  This means that we obtain an isomorphism
  $\HH^2(S_b)_\new^\tra \to \HH^2(X_b)^\tra$.
  In particular, the motive $\HH^2(S_b)_\new^\tra$ is abelian.
  To conclude that $\HH^2(S)$ is abelian, observe that
  $\HH^2(S) \cong \HH^2(S)_\new^\tra \oplus \HH^2(S)_\old^\tra \oplus \QQ(-1)^r$.
  The term $\HH^2(S)_\old^\tra$ is abelian,
  because it is the part coming from the Albanese surface,
  whose motive is abelian by definition.
 \end{proof}
\end{theorem}

\begin{corollary} 
 \label{T-MTC}
 Let $K$ be a finitely generated subfield of~$\CC$.
 Let $S$ be a surface of general type over~$K$
 with invariants $p_g = q = 2$ and dominant Albanese morphism.
 Assume that $S$ lies in one of the connected components
 of the Gieseker moduli space of surfaces of general type
 that contain a surface that is (semi-)isogenous to a product of curves.
 (That is, one of the types \textnumero~3,~4, 6, 8, or~9
 in \cref{state_of_the_art}.)
 Then the Tate and Mumford--Tate conjectures are true for~$S$.
 \begin{proof}
  We first prove the Mumford--Tate conjecture for~$S$.
  Let $A$ be the Albanese variety of~$S$.
  By \cref{motivated-k3-partner} there exists a K3~surface~$X$
  such that $\HH^2(S)_\new^\tra \cong \HH^2(X)^\tra$.
  Possibly after replacing~$K$ by a finitely generated extension
  we may assume that~$X$ is defined over~$K$.
  Hence the motive $\HH(S)$ is an object in the
  Tannakian subcategory of~$\Mot_K$ generated by~$\HH(A)$ and~$\HH(X)$.
  Therefore it suffices to prove the Mumford--Tate conjecture
  for $\HH(A) \oplus \HH(X)$.
  This follows from the main result of~\cite{Co16}.
  (See also \cite{Va08} and \cite{mtcaxa} for more general
  results on the Mumford--Tate conjecture for
  direct sums of abelian motives.)

  Recall that the Hodge conjecture is true for~$S$,
  by the Lefschetz-$(1,1)$ theorem.
  Therefore the Tate conjecture for~$S$ is true,
  since it follows from the conjunction of
  the Hodge conjecture and the Mumford--Tate conjecture.
  Indeed, if $\gamma \subset \Hl^2(S)(1)$ is a Tate class,
  then this means that it is fixed by $\Glc(S)$.
  We have just proven the Mumford--Tate conjecture for~$S$,
  so we know that $\Glc(S) \cong \GB(S) \otimes \QQl$.
  This means that $\gamma \in \Hl^2(S)(1) \cong \HB^2(S)(1) \otimes \QQl$
  is a $\QQl$-linear combination of
  $\GB(S)$-invariant classes in~$\HB^2(S)(1)$.
  Those $\GB(S)$-invariant classes are precisely Hodge classes,
  and by the Lefschetz-$(1,1)$ theorem
  we know that they are in the image of the cycle class map.
  We conclude that $\gamma$ is in the $\QQl$-span
  of the image of the $\ell$-adic cycle class map.
 \end{proof}
\end{corollary}

\printbibliography
\end{document}